\documentstyle[11pt]{article}

\input amssym.def
\input amssym

\begin{document}

\title{Torus actions on compact quotients} 
\author{Anton Deitmar\\ {\small Math. Inst. d. Univ., INF 288, 69126 Heidelberg, Germany}}
\date{}
\maketitle

\pagestyle{myheadings}
\markright{TORUS ACTIONS...}

\tableofcontents

\newfont{\spezial}{msbm10}

\def \a{{{\frak a}}}
\def \ad{{\rm ad}}
\def \al{\alpha}
\def \ar{{\alpha_r}}
\def \A{{\Bbb A}}
\def \Ad{{\rm Ad}}
\def \b{{{\frak b}}}
\def \bs{\backslash}
\def \B{{\cal B}}
\def \cent{{\rm cent}}
\def \C{{\bf C}}
\def \CA{{\cal A}}
\def \CB{{\cal B}}
\def \CE{{\cal E}}
\def \CF{{\cal F}}
\def \CG{{\cal G}}
\def \CH{{\cal H}}
\def \CM{{\cal M}}
\def \CN{{\cal N}}
\def \CP{{\cal P}}
\def \CQ{{\cal Q}}
\def \CO{{\cal O}}
\def \det{{\rm det}}
\def \e{\epsilon}
\def \fisch{\mbox{\spezial \char'156}}
\def \End{{\rm End}}
\def \Fx{{\frak x}}
\def \FX{{\frak X}}
\def \g{{{\frak g}}}
\def \ga{\gamma}
\def \Ga{\Gamma}
\def \h{{{\frak h}}}
\def \Hom{{\rm Hom}}
\def \ind{{\rm ind}}
\def \Im{{\rm Im}}
\def \Ind{{\rm Ind}}
\def \k{{{\frak k}}}
\def \K{{\cal K}}
\def \la{\lambda}
\def \lap{\triangle}
\def \m{{{\frak m}}}
\def \mod{{\rm mod}}
\def \Mat{{\rm Mat}}
\def \n{{{\frak n}}}
\def \name{\bf}
\def \N{\Bbb N}
\def \ord{{\rm ord}}
\def \O{{\cal O}}
\def \p{{{\frak p}}}
\def \ph{\varphi}
\def \prf{{\bf Proof: }}
\def \qed{\hfill {$\Box$} 

$ $

}
\def \Q{\Bbb Q}
\def \r{{\frak r}}
\def \res{{\rm res}}
\def \R{{\Bbb R}}
\def \Re{{\rm Re \hspace{1pt}}}
\def \ra{\rightarrow}
\def \rank{{\rm rank}}
\def \supp{{\rm supp}}
\def \SL{{\rm SL}}
\def \t{{{\frak t}}}
\def \T{{\Bbb T}}
\def \tr{{\hspace{1pt}\rm tr\hspace{1pt}}}
\def \vol{{\rm vol}}
\def \V{{\cal V}}
\def \z{\zeta}
\def \Z{\Bbb Z}

\newcommand{\rez}[1]{\frac{1}{#1}}
\newcommand{\der}[1]{\frac{\partial}{\partial #1}}
\newcommand{\binom}[2]{\left( \begin{array}{c}#1\\#2\end{array}\right)}

\newcounter{lemma}
\newcounter{corollary}
\newcounter{proposition}
\newcounter{theorem}

\newtheorem{conjecture}{\stepcounter{lemma} \stepcounter{corollary} 	
	\stepcounter{proposition}\stepcounter{theorem}Conjecture}[section]
\newtheorem{lemma}{\stepcounter{conjecture}\stepcounter{corollary}	
	\stepcounter{proposition}\stepcounter{theorem}Lemma}[section]
\newtheorem{corollary}{\stepcounter{conjecture}\stepcounter{lemma}
	\stepcounter{proposition}\stepcounter{theorem}Corollary}[section]
\newtheorem{proposition}{\stepcounter{conjecture}\stepcounter{lemma}
	\stepcounter{corollary}\stepcounter{theorem}Proposition}[section]
\newtheorem{theorem}{\stepcounter{conjecture} \stepcounter{lemma}
	\stepcounter{corollary}	\stepcounter{proposition}Theorem}[section]

$$ $$

 {\bf Introduction} 
 
Let $G$ denote a Lie group and $\Ga$ a uniform lattice in $G$.
We fix a maximal torus $T$ in $G$ and consider the action of $T$ on the compact 
quotient $\Ga \bs G$.
Assuming $T$ to be noncompact we 
will prove a Lefschetz formula relating compact orbits as local data to 
the action of the torus $T$ on a global cohomology theory (tangential 
cohomology).
Modulo homotopy, the compact orbits are parametrized by those conjugacy 
classes $[\ga]$ in $\Ga$ whose $G$-conjugacy classes meet $T$ in points 
which are 
regular in the split component.
Having a bijection between homotopy classes and conjugacy classes in the 
discrete group we will identify these two.
For a class $[\ga]$ let $X_\ga$ be the union of all compact 
orbits in that class. 
Then it is known that $X_\ga$ is a smooth 
submanifold and with $\chi_{_r}(X_\ga)$ we denote its de-twisted 
Euler characteristic (see sect. \ref{eulerchar}).
Note that $\chi_{_r}(X_\ga)$ is local, i.e. it can be expressed as 
the integral over $X_\ga$ of a canonical differential form 
(generalized Euler form).
On the other hand $\chi_{_r}(X_\ga)$ can be expressed 
as a simple linear combination of Betti numbers (see sect. \ref{eulerchar}).
Next, $\la_\ga$ will denote the volume of the orbit and $P_s$ the stable part of the Poincar\'e map around the orbit.
Then the number
$$
L (\ga) := \frac{\la_\ga \chi_{_r}(X_\ga)}{\det(1-P_s)}
$$
will be called the Lefschetz number of $[\ga]$ (compare \cite{gui}).
The class $[\ga]$ defines a point $a_\ga$ in the split part $A$ of the torus $T$ 
modulo the action of the Weyl group.
In the case when the Weyl group has maximal size (for example when $T$ 
is maximally split) our Lefschetz formula is an equality of distributions:
$$
\sum_{[\ga]} L(\ga) \delta_{a_\ga} = \tr (. | H^*(\CF)),
$$
where $H^*$ is the tangential cohomology of the unstable/neutral foliation $\CF$ induced by the torus action.
In \cite{DG} a similar formula is proven to hold up to a smooth function in the case of a flow.
The present paper extends results of Andreas Juhl \cite{Ju}, \cite{Schub} 
in the real rank one case.
See also \cite{Osb1}, \cite{Osb2}.

$$ $$

\section{Euler-Poincar\'e functions} 
In this section and the next we list some technical results for the convenience of the reader. Let $G$ denote a real reductive group of inner type \cite{Wall-rg1} and fix a maximal compact subgroup $K$.
Let $(\tau ,V_\tau)$ be a finite dimensional unitary representation of $K$ and write $(\breve{\tau},V_{\breve{\tau}})$ for the dual representation. 
Assume that $G$ has a compact Cartan subgroup $T \subset K$.
Let $\g_0 = \k_0 \oplus \p_0$ be the polar decomposition of the real Lie algebra $\g_0$ of $G$ and write $\g = \k +\p$ for its complexification.
Choose an ordering of the roots $\Phi(\g ,\t)$ of the pair $(\g ,\t)$. This choice induces a decomposition $\p = \p_- \oplus \p_+$.

\begin{proposition} \label{existf}
For $(\tau ,V_\tau)$ a finite dimensional representation of $K$ there is a compactly supported smooth function $f_\tau$ on $G$ such that for every irreducible unitary representation $(\pi ,V_\pi)$ of $G$ it holds:
$$
\tr\ \pi (f_\tau) = \sum_{p=0}^{\dim (\p)} (-1)^p \dim (V_\pi \otimes \wedge^p\p \otimes V_{\breve{\tau}})^K.
$$
\end{proposition}
\prf
\cite{onsome}.
\qed

\begin{proposition} \label{orbitalint}
Let $g$ be a semisimple element of the group $G$. If $g$ is not elliptic, then the orbital integral $\O_g(f_\tau)$ vanishes. If $g$ is elliptic we may assume $g\in T$, where $T$ is a Cartan in $K$ and then we have
$$
\O_g(f_\tau) = \frac{{\tr\ \tau(g)}|W(\t ,\g_g)| \prod_{\alpha \in \Phi_g^+}(\rho_g ,\alpha)}{[G_g:G_g^0]c_g},
$$
where $c_g$ is Harish-Chandra's constant, it does only depend on the centralizer $G_g$ of g. Its value is given for example in \cite{hitors}.
\end{proposition}

\prf
\cite{onsome}.
\qed

\begin{proposition}
For the function $f_{{\sigma}}$ we have for any $\pi \in \hat{G}$:
$$
\tr \ \pi(f_{{\sigma}}) = \sum_{p=0}^{\dim \ \g /\k}(-1)^p \dim \ 
{\rm Ext}_{(\g ,K)}^p (V_\sigma ,V_\pi),
$$
i.e. $f_{{\sigma}}$ gives the Euler-Poincar\'e numbers of the $(\g ,K)$-modules $(V_\sigma ,V_\pi)$, this justifies the name Euler-Poincar\'e function.
\end{proposition}

\prf
\cite{onsome}.
\qed

\section{De-twisted Euler characteristics}\label{eulerchar}
Let ${\cal C}^+$ denote the category of complexes of $\C$-vector spaces which are zero in negative indices and have degreewise finite dimensional cohomology, i.e.
the dimension of $H^j(E)$ is finite for all $j$. Let $\K^+$ denote the weak Grothendieck group of ${\cal C}^+$, i.e. $\K^+$ is the abelian group generated by all isomorphism classes of objects modulo the relations $A=B+C$, whenever any object in $A$ is isomorphic to the direct sum of an object in $B$ and one in $C$.
An element $E=E_+-E_-$ of $\K^+$ is called a {\bf virtual complex}. 
Define the {\bf de-twist} of an element $E$ of $\K^+$ as
$$
E' = \sum_{k=0}^\infty E[-k],
$$
where $E[k]_j=E_{k+j}$. Since the sum is degreewise finite this defines a new element of $\K^+$. The higher de-twists are defined inductively, so $E^{(0)}=E$ and $E^{(r+1)}={E^{(r)}}'$.

We need to extend   the notion of an {\bf Euler characteristic} to infinite virtual complexes by
$$
\chi (E) = \sum_{k=0}^\infty (-1)^k \dim H^k(E),
$$
provided $\dim H^k(E) = \dim H^k(E_+) - \dim H^k(E_-)$ vanishes for almost all $k$.

Call a virtual complex {\bf cohomologically finite} if $ \dim H^j (E) = 0 $ for large j, in other words, the total cohomology H(E) is finite dimensional.

{\bf Observation:} {\it Let the virtual complex E be cohomologically finite and assume that the Euler characteristic $ \chi (E) $ vanishes. Then the de-twist $E'$ is cohomologically finite}.

So start with a cohomologically finite virtual complex E. If
$ E^{(1)}, \dots , E^{(r)} $ are cohomologically finite we have
$$
\chi ( E^{(0)}) = \dots = \chi (E^{(r-1)}) = 0
$$
and
$$
\chi (E^{(r)}) = (-1)^r \sum_{j=0}^\infty \binom{j}{r}
(-1)^j dim H^j(E).
$$
This is easily proven by induction on r. This motivates the following Definition: The \label{higher Euler characteristic} \index{higher Euler characteristic}{\bf r-th de-twisted Euler characteristic} of a cohomologically finite virtual complex E is defined by
$$
\chi_r (E) := (-1)^r \sum_{j=0}^\infty \binom{j}{r} (-1)^j dim H^j(E).
$$ \index{$\chi_r (E)$}

To every compact manifold M we now can attach a sequence of Euler numbers

$$
\chi_{_0}(M), \dots ,\chi_{_n}(M),
$$
where $n$ is the dimension of M. 
The most significant of these is, as we shall see, the first 
nonvanishing one, so define the 
\index{generic Euler number}{\bf generic Euler number} of M as
$$
\chi_{_{gen}}(M) = \chi_{_r}(M), \ \ \ 
{\rm where}\ r\ {\rm is\ the\ least\ index\ with}\
                 \chi_{_r}(M) \neq 0.
$$
\begin{proposition} Let M,N be compact manifolds. We have
$$
\chi_{_{gen}}(M\times N) = \chi_{_{gen}}(M) \chi_{_{gen}}(N).
$$
\end{proposition}

\prf
See \cite{D-Prod}.
\qed

To give another example of a situation in which higher Euler characteristics occur we will describe a situation in Lie algebra cohomology which will show up later.

We consider a short exact sequence
$$
o \ra \n \ra \l \ra \a \ra 0
$$
of finite dimensional complex  Lie algebras where $\a$ is abelian.
In such a situation a $\l$-module $V$ is called 
{\bf acceptable} \index{acceptable $\l$-module}, 
if the $\a$-module $H^q(\n ,V)$ is finite dimensional. 
Note that $V$ itself needn't be finite dimensional.

{\bf Example 1.:} Any finite dimensional $\l$-module will be acceptable.

{\bf Example 2.:} Let $\g_0$ denote the Lie algebra of a semisimple Lie group $G$ of the Harish-Chandra class, i.e. $G$ is connected and has a finite center. Let $K$ be a maximal compact subgroup of $G$ and let $G=KAN$ be an Iwasawa decomposition of $G$. Write the corresponding decomposition of the complexified Lie algebra as $\g = \k \oplus \a \oplus \n$. Now let $\l = \a \oplus \n$ with the structure of a subalgebra of $\g$. Consider an admissible $(\g ,K)$-module $V$. A theorem of [HeSchm] assures us that $V$ then is an acceptable $\l$-module.

\begin{proposition} Let
$$
o \ra \n \ra \l \ra \a \ra 0
$$
be an exact sequence of finite dimensional complex Lie algebras. Assume that the Lie algebra $\a$ is abelian. Let $V$ be an acceptable $\l$-module then with $r=\dim(\a)$ we have
$$
\chi_{_0}(H^*(\l ,V)) = \dots = \chi_{_{r-1}}(H^*(\l ,V)) = 0,
$$
and
$$
\chi_{_r}(H^*(\l ,V)) = \chi_{_0}(H^*(\n ,V)^\a),
$$
where $H^*(\n ,V)^\a$ denotes the $\a$-invariants in $H^*(\n ,V)$.
\end{proposition}

\prf
\cite{onsome}.
\qed

\section{The Lefschetz formula}
Let $G$ be a connected Lie group and $\Ga \subset G$ a uniform lattice.
Note that the existence of $\Ga$ forces $G$ to be unimodular.
Fix a Haar measure on $G$ and
consider the representation of $G$ on the Hilbert space $L^2(\Ga \bs G)$ 
given by $R(g) \ph(x)=\ph (xg)$.
For any smooth compactly supported function $f$ on $G$ define 
$R(f)\ph(x):=\int_Gf(y)\ph(xy)dy$, then a calculation shows that 
$R(f)$ is an integral operator with smooth kernel 
$k(x,y) = \sum_{\ga\in\Ga}f(x^{-1}\ga y)$.
From this it follows that $R(f)$ is a trace class operator.
Since this holds for any $f$ we conclude that $L^2(\Ga \bs G)$ 
decomposes under $G$ as a discrete sum of irreducibles with finite 
multiplicities:
$$
L^2(\Ga \bs G) = \bigoplus_{\pi\in\hat{G}}N_\Ga(\pi) \pi .
$$

It follows that $\tr R(f) = \sum_{\pi\in\hat{G}}N_\Ga(\pi) \tr \pi(f)$.
On the other hand, the trace of $R(f)$ equals the integral over the 
diagonal of the kernel, so
\begin{eqnarray*}
\tr R(f) &=& \int_{\Ga \bs G} k(x,x) dx\\
	&=& \sum_{[\ga]}\vol(\Ga_\ga \bs G_\ga) \CO_\ga(f),
\end{eqnarray*}
where $\CO_\ga(f) := \int_{G_\ga \bs G}f(x^{-1}\ga x)dx$ is the 
orbital integral.
Note that this expression depends on the choice of a Haar measure 
on $G_\ga$.
So we state the Selberg trace formula as
$$
\sum_{\pi\in\hat{G}}N_\Ga(\pi) \tr \pi(f)= 
\sum_{[\ga]}\vol(\Ga_\ga \bs G_\ga) \CO_\ga(f).
$$

From now on we will {\bf assume}:
\begin{itemize}
\item[{\bf (A1)}]
$G$ is a direct product:
$$
G\cong H\fisch R
$$
of an abelian Lie group $R$ and a semisimple connected Lie group $H$ with finite center.
\end{itemize}

For the following fix a maximal torus $T$ of $H$, write $T=AB$, 
where $A$ is the split component and $B$ is compact.
Let $P=MAN$ a parabolic then $B\subset M$.
Let $A^{reg}$ be the set of regular elements of the split torus $A$.
Since $H$ acts on $R$ it acts on the unitary dual $\hat{R}$.
Our second assumption is
\begin{itemize}
\item[{\bf (A2)}]\label{A3}
Any element of $A^{reg}M$ acts freely on 
$R - \{ 0\}$ and on $\hat{R} - \{ triv\}$.
\end{itemize}

For any $\tau\in\hat{R}$ let $H_\tau$ be its stabilizer in $H$. 
For the trivial representation we clearly have
$H_{triv} =H$.
The condition (A2) says that for any nontrivial $\tau\in\hat{R}$ we 
have $H_\tau \cap A^{reg}M=\emptyset$.

$ $

{\bf Example:}
Clearly any semisimple connected $G$ with finite center would give an example 
but there are also a lot on nonreductive examples such as the following:
Let $R := \Mat_2(\R)$ with the addition, $H:=\SL_2(\R)$ and 
let $H$ act on $R$ by matrix multiplication from the left. 
Let $G:=H\fisch N$ and 
$T:=\left\{\left( \begin{array}{cc}a&{}\\{}&a^{-1}\end{array}\right)\right\}$.
It is easily seen that our assumptions are satisfied in this case.

$ $

We will only consider uniform lattices of the form 
$\Ga =\Ga_H\fisch \Ga_R$, where $\Ga_H$ and $\Ga_R$ are 
uniform lattices in $H$ and $R$.
We will further assume $\Ga_H$ to be {\bf weakly neat},
this means, $\Ga_H$ is a cocompact torsion free discrete 
subgroup of $H$ which is such that for any $\ga\in\Ga_H$ 
the adjoint $\Ad(\ga)$, acting on the Lie algebra of $H$ 
does not have a root of unity $\neq 1$ as an eigenvalue.
Any arithmetic group has a weakly neat subgroup of finite 
index \cite{Bor}.

$ $

{\bf Example:}
Take up the above example and let $D$ denote a quaternion division 
algebra over $\Q$ which splits over $\R$.
So we have $D\hookrightarrow GL_2(\R)$ and $D^1\hookrightarrow SL_2(\R)$, 
where $D^1$ is the set of elements of reduced norm $1$.
Let $\CO$ denote an order in $D$ and $\CO^1 :=\CO\cap D^1$.
Then $\Ga :=\CO^1 \fisch \CO$ is a uniform lattice in 
$SL_2(\R)\fisch\Mat_2(\R)$.

$ $

Write the real Lie algebras of $G, H, M, A, N, R$ as 
$\g_0, \h_0, \m_0, \a_0, \n_0, \r_0$ and their complexifications 
as $\g, \h, \m, \a, \n, \r$.
Let $\Phi (\h ,\a)$ denote the set of roots of the pair 
$(\h ,\a)$.
The choice of the parabolic $P$ amounts to the same as a 
choice of a set of positive roots $\Phi^+ (\h ,\a)$.
Let $A^-\subset A$ denote the negative Weyl chamber corresponding 
to that ordering, i.e. $A^-$ consists of  all $a\in A$ which act 
contractingly on the Lie algebra $\n$.
Further let $\overline{A^-}$ be the closure of $A^-$ in $G$, 
this is a manifold with boundary.
Let $K_M$ be a maximal compact subgroup of $M$. We may suppose 
that $K_M$ contains $B$.
Fix an irreducible unitary representation $(\tau ,V_\tau)$ of $K_M$.
Let $K$ be a maximal compact subgroup of $H$. We may assume  $K\supset K_M$.

Since $\Ga_H$ is the fundamental group of the Riemannian manifold
$$
X_{\Ga_H} = \Ga_H \bs X = \Ga_H \bs H/K
$$
 it follows that we have a 
canonical bijection of the homotopy classes of loops:
$$
[S^1 : X_{\Ga_H} ] \ra \Ga_H / {\rm conjugacy}.
$$
For a given class $[\ga]$ let $X_\ga$ denote the union of all 
closed geodesics in the corresponding class in $[S^1 : X_\Ga ]$. 
Then $X_\ga$ is a smooth submanifold of $X_{\Ga_H}$ \cite{DKV}.
Let $\chi_{_r}(X_\ga)$ denote the $r$-fold de-twisted 
Euler characteristic of $X_\ga$, where $r=\dim A$.

Let $\CE_P(\Ga)$ denote the set of all conjugacy classes 
$[\ga]$ in $\Ga$ such that $\ga_H$ is in $H$ conjugate to an 
element $a_\ga b_\ga$ of $A^- B$.

Take a class $[\ga]$ in $\CE_P(\Ga)$. Modulo conjugation assume 
$\ga\in T=AB$, then the centralizer $\Ga_{H,\ga}$ projects to a 
lattice $\Ga_{A,\ga}$in the split part $A$.
Let $\la_\ga$ be the covolume of this lattice.
Normalize the measure on $R$ such that $\vol(\Ga_R\bs R)=1$.

\begin{theorem} \label{ersteversion} (Lefschetz formula, first version)
Let $\ph$ be compactly supported on $\overline{A^-}$, $\dim G$-times 
continuously differentiable and suppose $\ph$ vanishes on the boundary 
to order $\dim G +1$.
Then we have that the expression
$$
\sum_{\begin{array}{c}\pi \in \hat{G}\\ \pi |_R \equiv 1\end{array}} N_\Ga (\pi) 
\sum_{p,q}(-1)^{p+q}
\int_{A^-} \ph (a) \tr(a|(H^q({\n},\pi )\otimes \wedge^p\p_M\otimes V_{\breve{\tau}})^{K_M}) da
$$
equals
$$
 (-1)^{\dim(N)} \sum_{[\ga]\in {\cal E}_P(\Ga)}\la_\ga \chi_{_r}(X_\ga) \frac{\ph(a_\ga)\tr \tau(b_\ga)}{\det(1- a_\ga b_\ga |{\n})\det(1-\ga | \r)}.
$$
\end{theorem}

\prf
Let $H$ act on itself by conjugation, write 
$h.x = hxh^{-1}$, write $H.x$ for the orbit, so 
$H.x = \{ hxh^{-1} | h\in H \}$ as well as 
$H.S = \{ hsh^{-1} | s\in S , h\in H \}$ for any subset $S$ of $H$.
We are going to consider functions that are supported on the closure 
of the set $H.(MA^-)$.
At first let $f_\tau$ be the Euler-Poincar\'e function defined on 
$M$ attached to the representation $(\tau ,V_\tau)$ of $K_M$.
Next fix a smooth function $\eta$ on $N$ which has compact support, 
is positive, invariant under $K_M$ and satisfies $\int_N\eta(n) dn =1$.
Given these data let $\phi = \phi_{\eta ,\tau ,\ph} : H\ra \C$ 
be defined by
$$
\phi (kn ma (kn)^{-1}) := \eta (n) f_\tau(m) \frac{\ph(a)}{\det(1-(ma)|\n)},
$$
for $k\in K, n\in N, m\in M, a\in\overline{A^-}$.
Further $\phi(h)=0$ if $h$ is not in $H.(M\overline{A^-})$.

Next choose any compactly supported positive function $\psi$ on $R$ with 
$\int \psi =1$.
Let $\Phi (h,r) := \phi(h)\psi(r)$.
We will plug $\Phi$ into the trace formula.
For the geometric side let $\ga =(\ga_H,\ga_R)\in \Ga$. We have to calculate the orbital integral:
$$
\CO_\ga (\Phi) = \int_{G_\ga \bs G} \Phi(x^{-1}\ga x) dx.
$$

Now let $x=(h,r)\in G$ and compute
$$
x^{-1}\ga x = (h^{-1}\ga_H h, r+h^{-1}\ga_R -h^{-1}\ga_H hr).
$$
So (h,r) lies in the centralizer $G_\ga$ iff $h\in H_{\ga_H}$ and $r\in R$ satisfies
$$
(1-h^{-1})\ga_R = (1-\ga_H^{-1})r.
$$
Note that by (A2) to any $\ga$ such that $\ga_H$ is conjugate to an element of $A^{reg}M$, and to any $h\in H_{\ga_H}$ such an $r$ exists and is unique.
But this condition on $\ga$ is satisfied if $\ph(h^{-1}\ga_H h)\neq 0$.
So suppose $\ga_H$ is in $H.(A^{reg}M)$.
In this case we have the integration rule
$$
\int_{G_\ga \bs G} f(g) dg = \int_{H_{\ga_H}\bs H} \int_R f(h,r) dr dh.
$$
This is proven by showing that the right hand side is in fact $G$-invariant.
We compute
$$
\int_R \Phi((h,r)^{-1}\ga(h,r))dr = \frac{\ph(h^{-1}\ga_H h)}{\det(1-\ga_H |\r)},
$$
from which we see that the geometric side of the trace formula coincides with our claim.

Now for the spectral side let $\pi\in\hat{G}$ then the restriction of $\pi$ to $R$ is a direct integral over $\hat{R}$.
The irreducibility of $\pi$ implies that the corresponding measure is supported on a single orbit $o$ of the $H$-action on $\hat{R}$.
So we have
$$
\pi |_R = \int_o V_\pi (\tau)dm(\tau),
$$
where $m$ is a scalar valued measure and $V_\pi(\tau)$ is a multiple of $\tau$.
Fix $\tau_0\in o$ then the stabilizer $H_{\tau_0}$ acts trivially on $\tau_0$ and not only its class since by $\dim \tau_0 =1$ these two notions coincide.
It follows that as $H_{\tau_0}$-modules we have $V_\pi(\tau) \cong \eta \otimes \tau$ for some representation $\eta$ of $H_{\tau_0}$.
The measure $m$ induces a measure on $G_{\tau_0}\bs G$ also denoted $m$ which is quasi-invariant.
It follows that $\pi = \ind_{H_{\tau_0}\fisch R}^G(\eta\otimes\tau)$ and hence $\eta$ must be irreducible since $\pi$ is.
Let $\la (x,y)$ denote the Radon-Nikodym derivative of the translate $m_x$ with respect to $m$.
We conclude that $\pi(\Phi)$ is given as an integral operator on $G_{\tau_0}\bs G$ with kernel
$$
k(x,y)= \int_{G_{\tau_0}} \Phi(x^{-1}zy)\la(x^{-1}zy,x)^{\rez{2}} (\eta \otimes \tau_o)(z) dz.
$$
From this we get
$$
\tr \pi(\Phi) = \int_{G_{\tau_0}} \tr(\int_{G_{\tau_0}}\Phi(x^{-1}zx)\la(x^{-1}zx,x)^{\rez{2}}(\eta\otimes\tau_0)(z)dz)dx.
$$
Consider the term $\Phi(x^{-1}zx) = \ph(x_H^{-1}z_Hx_H)\psi(\dots)$.
By (A2) this expression vanishes unless $\tau_0$ is the trivial character.
In the case $\tau_0=triv$ it follows that $\pi(R)=1$, so $\pi$ may be viewed as an element of $\hat{H}$.

To evaluate $\tr \pi(\Phi)$ further we will employ the Hecht-Schmid character formula \cite{HeSch}.
For this let
$$
(MA)^- = {\rm interior\ in\ } MA {\rm \ of\ the\ set}
$$ $$
\left\{ g\in MA | \det (1-ga | \n) \geq 0 {\rm for\ all\ } \a\in A^- \right\}.
$$
The character $\Theta_\pi^G$ of $\pi\in\hat{G}$ is a locally integrable function on $G$.
In \cite{HeSch} it is shown that for any $\pi\in\hat{H}$, denoting by $\pi^0$ the underlying Harish-Chandra module we have that all Lie algebra cohomology groups $H^p(\n ,\pi^0)$ are Harish-Chandra modules for $MA$.
The main result of \cite{HeSch} is that for $ma\in (MA)^- \cap H^{reg}$, the regular set, we have
$$
\Theta_\pi^H(ma) = \frac{\sum_{p=0}^{\dim \n}(-1)^p \Theta_{H_p(\n ,\pi^0)}^{MA}(ma)}{\det (1-ma |\n)}.
$$
Let $f$ be supported on $H.(MA^-)$, then the Weyl integration formula states that
$$
\int_H f(x) dx = \int_{H/MA} \int_{MA^-} f(hmah^{-1}) |det(1-ma|\n\oplus\bar{n})dadmdh.
$$
So that for $\pi\in\hat{H}$:
\begin{eqnarray*}
\tr \pi(\phi) &=& \int_H \Theta_\pi^H (x) \phi (x) dx\\
	&=& \int_{MA^-} \Theta_\pi^H (ma) f_\tau(m)\ph(a) |\det(1-ma|\bar{\n})|dadm\\
	&=& (-1)^{\dim N} \int_{MA^-} f_\tau (m) \Theta_{H^*(\n ,\pi^0)}^{MA}(ma)\ph(a) da dm,
\end{eqnarray*}
where we have used the isomorphism $H_p(\n ,\pi^0)\cong H^{\dim N -p}(\n ,\pi^o)\otimes \wedge^{top}\n$.
This gives the claim.
\qed

In the second version of the Lefschetz formula we want to substitute the character of the representation $\tau$ by an arbitrary central function on $K_M$.
A smooth function $f$ on $K_M$ is called {\bf central} if $f(kk_1k^{-1})=f(k_1)$ for all $k,k_1\in K_M$.
Since $B$ is a Cartan subgroup of the compact group $K_M$, any $k\in K_M$ is conjugate to some element of $B$ so the restriction gives in isomorphism from the space of smooth central functions on $K_M$ to the space of smooth functions on $B$, invariant under the Weyl group.
Hence we are led to consider Weyl group invariant functions on $T$.

Let $\CA$ denote the convolution algebra of all $W(H,T)$-invariant smooth functions on $T$ with compact support.
Let $S\subset T$ be the set of all $ab$ with singular $a$-part.

For any $t=ab$ in $T$ let $\n_t$ be the space of all $X\in \ad (t)\g$ on which $t$ acts contractingly.
Then $\n_t$ is a nilpotent Lie subalgebra of $\g$.

\begin{theorem} (Lefschetz formula, second version)
Let $\ph\in\CA$ and suppose $\ph$ vanishes on the singular set to order $\dim G +1$ then the expression
$$
\sum_{\begin{array}{c}\pi \in \hat{G}\\ \pi |_R \equiv 1\end{array}} N_\Ga (\pi) 
\sum_q (-1)^q
\int_{T/W(H,T)} \ph (t) \tr(t|H^q({\n_t},\pi )) dt
$$
equals
$$
 (-1)^{\dim(N)} \sum_{[\ga]\in {\cal E}_P(\Ga)}\la_\ga \chi_{_r}(X_\ga) \frac{\ph(t_\ga)}{\det(1- t_\ga |\p_M\oplus\n_{h_\ga})\det(1-\ga |\r)}.
$$
\end{theorem}

\prf
Extend $b\mapsto \ph(ab)$ to a central function on $K_M$. Then expand $\ph$ 
into $K_M$-types:
$$
\ph (ab) = \sum_{\tau\in\hat{K_M}} c_\tau \tr \tau(b) \ph_\tau(a),
$$
since $\ph$ is smooth the coefficients $c_\tau$ are rapidly decreasing 
so the expressions of Theorem \ref{ersteversion} when plugging 
in $\ph_\tau |_{A^-}$ converge to 
$$
\sum_{\begin{array}{c}\pi \in \hat{G}\\ \pi |_R \equiv 1\end{array}} N_\Ga (\pi)
\sum_{p,q} (-1)^{p+q}
 \int_{T/W(H,T)} \ph (t) \tr(t|H^q({\n_t},\pi )\otimes \wedge^p\p_M) dt, 
$$
which equals
$$
 (-1)^{\dim(N)} \sum_{[\ga]\in {\cal E}_H(\Ga)}\la_\ga \chi_{_r}(X_\ga) \frac{\ph(t_\ga)}{\det(1- t_\ga |\n_{h_\ga})\det(1-\ga |\r)}.
$$
Now replace $\ph(t)$ by $\ph(t)/\det(1-t|\p_M)$ which gives the claim.
\qed

At last we also mention a reformulation in terms of relative Lie algebra cohomology.
Again, fix a parabolic $P=MAN$ and now fix also a finite dimensional irreducible representation $(\sigma ,V_\sigma)$ of $M$. 

\begin{theorem} \label{dritteversion} (Lefschetz formula, third version) 
Let $\ph$ be compactly supported on $\overline{A^-}$, $\dim G$-times 
continuously differentiable and suppose $\ph$ vanishes on the boundary 
to order $\dim G +1$.
Then we have that the expression
$$
\sum_{\begin{array}{c}\pi \in \hat{G}\\ \pi |_R \equiv 1\end{array}} N_\Ga (\pi)
\sum_q(-1)^q
 \int_{A^-} \ph (a) \tr(a|H^q(\m\oplus{\n},K_M,\pi \otimes V_{\breve{\sigma}}) 
$$
equals
$$
 (-1)^{\dim(N)} \sum_{[\ga]\in {\cal E}_P(\Ga)}\la_\ga \chi_{_r}(X_\ga) \frac{\ph(a_\ga)\tr \sigma(b_\ga)}{\det(1- a_\ga b_\ga |{\n})\det(1-\ga |\r)}.
$$
\end{theorem}

\prf
Extend $V_{\breve{\sigma}}$ to a $\m \oplus {\n}$-module by letting ${\n}$ act trivially. We then get
$$
H^p({\n},\pi^0) \otimes V_{\breve{\sigma}} \cong H^p({\n},\pi^0 \otimes V_{\breve{\sigma}}).
$$

The $(\m ,K_M)$-cohomology of the module $H^p({\n},\pi^0 \otimes V_{\breve{\sigma}})$ is the cohomology of the complex $(C^*)$ with
\begin{eqnarray*}
C^q &=& {\rm Hom}_{K_M}(\wedge^q\p_M ,H^p({\n},\pi^0)\otimes V_{\breve{\sigma}})
\\
        &=& (\wedge^q\p_M \otimes H^p({\n},\pi^0)\otimes V_{\breve{\sigma}})^{K_M},
\end{eqnarray*}
since $\wedge^p\p_M$ is a self-dual $K_M$-module. Therefore we have an isomorphism of virtual $A$-modules:
$$
\sum_q (-1)^q (H^p({\n},\pi^0)\otimes \wedge^q\p_M \otimes V_{\breve{\sigma}} )^{K_M}
\cong
\sum_q (-1)^q H^q (\m ,K_M,H^p({\n},\pi^0 \otimes V_{\breve{\sigma}})).
$$

Now one considers the Hochschild-Serre spectral sequence in the relative case for the exact sequence of Lie algebras
$$
0 \ra {\n}\ra \m \oplus {\n}\ra \m \ra 0
$$
and the $(\m\oplus {\n},K_M)$-module $\pi \otimes V_{\breve{\sigma}}$. We have
$$
E_2^{p,q} = H^q(\m ,K_M ,H^p({\n},\pi^0\otimes V_{\breve{\sigma}}))
$$
and
$$
E_\infty^{p,q} = {\rm Gr}^q(H^{p+q}(\m \oplus {\n},K_M ,\pi^0 \otimes V_{\breve{\sigma}})).
$$
Now the module in question is just
$$
\chi(E_2) = \sum_{p,q} (-1)^{p+q} E_2^{p,q}.
$$
Since the differentials in the spectral sequence are $A$-homomorphisms this equals $\chi(E_\infty)$.
So we get an $A$-module isomorphism of virtual $A$-modules
$$
\sum_{p,q} (-1)^{p+q} (H^p({\n},\pi^0)\otimes \wedge^q\p_M \otimes V_{\breve{\sigma}})^{K_M} \cong \sum_j (-1)^j H^j (\m\oplus {\n},K_M, \pi^0 \otimes V_{\breve{\sigma}}).
$$
The claim follows.
\qed

\section{Geometric interpretation}
Now consider the first version of the Lefschetz formula in the case $R=0$.
The representation $\tau$ defines a homogeneous vector bundle 
$E_\tau$ over $G/K_M$ and by homogeneity 
this pushes down to a locally homogeneous bundle over $\Ga \bs G/K_M = {}_MX_\Ga$.
The tangent bundle $T(_MX_\Ga)$ can be described in this way as stemming from the representation of $K_M$ on
$$
\g /\k_M \cong \a \oplus \p_M \oplus \n \oplus\bar{\n}.
$$
We get a splitting into subbundles
$$
T(_MX_\Ga) = T_c \oplus T_n \oplus T_u \oplus T_s.
$$
These bundles can be characterized by dynamical properties: the action of $A\cong\R^r$ is furnished with a positive time direction given by the positive Weyl chamber $A^+$.
Then $T_s$, the {\bf stable part} is characterized by the fact that $A^+$ acts contractingly on $T_s$.
On the {\bf unstable part} $T_u$ the opposite chamber $A^-$ acts contractingly.
$T_c$, the {\bf central part} is spanned by the "flow" $A$ itself and $T_n$ is an additive {\bf neutral part}.
Note that $T_n$ vanishes if we choose $H$ to be the maximal split torus.
The bundle $T_n\oplus T_u$ is integrable, so it defines a foliation $\CF$.
To this foliation we have the tangential cohomology $H^*(\CF)$ and also for its $\tau$-twist: $H^*(\CF\otimes \tau)$.
The flow $A$ acts on the tangential cohomology whose alternating sum we will consider as a virtual $A$-module. 
For any $\ph\in C_c^\infty(\overline{A^-})$ we define $L_\ph = \int_{A^-}\ph(a) (a|H^*(\CF\otimes\tau)) da$ as a virtual operator on $H^*(\CF\otimes\tau)$.
Then we have

\begin{proposition}
Under the assumptions of theorem \ref{ersteversion} the virtual operator $L_\ph$ is of trace class and the RHS of Theorem \ref{ersteversion}
 can be written as
$$
\sum_q(-1)^q \tr (L_\ph |H^q(\CF\otimes\tau)).
$$
\end{proposition}

{\small Math. Inst. d. Univ., INF 288, 69126 Heidelberg, Germany}

\end{document}